\documentclass[twoside]{article}

\usepackage{amsmath}
\usepackage{amsthm}
\usepackage{amssymb}
\usepackage{latexsym}
\usepackage{array}

\usepackage{times} 
\usepackage{graphicx} 
\usepackage{color}
\usepackage{mysects}  
\usepackage{fancyhdr} 

\newtheorem{thm}{\bf Theorem\bf}[section]
\newtheorem{prop}[thm]{\bf Proposition\bf}
\newtheorem{cor}[thm]{\bf Corollary\bf}
\newtheorem{lem}[thm]{\bf Lemma\bf}

\newtheorem{rem}[thm]{\bf Remark\bf}

\def\proof{{\noindent\bf Proof. }}

\begin{document}

\setcounter{section}{0}
\setcounter{subsection}{0}
\setcounter{footnote}{0}
\setcounter{page}{1}

\thispagestyle{plain}
\pagestyle{fancy}
\fancyhead{}
\fancyhead[LE]{\thepage}
\fancyhead[RO]{\thepage}
\fancyhead[CO]{Real zeros of holomorphic Hecke cusp forms\hfill}
\fancyhead[CE]{\hfill Amit Ghosh and Peter Sarnak}
\fancyfoot[C]{}
\ \vskip 1cm
\renewcommand{\thefootnote}{ } 
\renewcommand{\footnoterule}{{\hrule}\vspace{3.5pt}} 

\vskip -1.0cm
{\Large\bf\noindent Real zeros of holomorphic Hecke cusp forms }
\vskip 0.2in
{\noindent \large Amit Ghosh and Peter Sarnak}

\vskip 0.5cm
{
\centerline{
\hbox{
\begin{minipage}{4.0in}
{\small  {\bf Abstract.}\quad This note is concerned with the zeros of holomorphic Hecke cusp forms of large weight on the modular surface. The zeros of such forms are symmetric about three geodesic segments and we call those zeros that lie on these segments, real. Our main results give estimates for the number of real zeros as the weight goes to infinity.
}
\end{minipage}
}  
}  
}
\vskip 0.5cm
\hrule
\footnote{Version of \today}
\renewcommand{\thefootnote}{\arabic{footnote}\quad } 
\setcounter{footnote}{0} 

\setcounter{section}{0}\setcounter{equation}{0}
\vskip 0.1in

\section{Introduction.}\label{Intro} 


For $k$ an even integer, the space of holomorphic forms of weight $k$ for the full modular group $\Gamma = PSL(2,\mathbb{Z})$ is of dimension $\frac{k}{12} + O(1)$ (see \cite{Ser} for exact definitions as well as other basic facts). Such a form $F$ has $\frac{k}{12} + O(1)$ zeros in $\mathfrak{X} = \Gamma\backslash\mathcal{H}$. More precisely,
\begin{equation}
\nu_{\infty}(F) + \frac{\nu_{i}(F)}{2} + \frac{\nu_{\rho}(F)}{3} + \sum_{p \in \mathfrak{X}}  \nu_{p}(F) = \frac{k}{2}
\end{equation}
where $\nu_{p}(F)$ is the order of vanishing of $F$ at the point $p$. Here $\infty$, $i$ and $\rho$ are points in $\mathfrak{X}$ depicted in the familiar Figure 1. Other than vanishing at $z=i$ and $z=\rho$ when forced by (1) for $k$ in various progressions modulo 12, the distribution of the zeros of such an $F$ is not restricted. However, for the arithmetically interesting case of $F$ being a Hecke eigenform, which we will assume henceforth, there are constraints on the distribution of the zeros. In particular, if $F$ is an Eisenstein series $E_{k}(z)$, then it has been shown by Rankin and Swinnerton-Dyer \cite{RanSD} that all of its zeros are on the geodesic segment $\delta_{3}$ in Figure 1 (there have been many generalizations of this result to functions constructed from Eisenstein series (see \cite{DukeJenkins}) and to other Fuchsian groups (see \cite{Hahn})). For the rest of the Hecke eigenforms, namely the cusp forms, which we denote by $f$, the distribution of the zero set $\mathcal{Z}(f)$ (counted with multiplicities) is very different. By definition such a $f$ vanishes at the cusp $z = i\infty$ and from the well known properties of the Hecke operators it follows that $\nu_{\infty}(f)=1$. The expansion of $f$ at $i\infty$ takes the form
\begin{equation}
f(z) = \sum_{n=1}^{\infty} a_{f}(n) e(nz),
\end{equation}
where we normalize $f$ with $a_{f}(1)=1$ and write $e(x) = exp(2\pi ix)$.

\centerline{\includegraphics[width=3.5in]{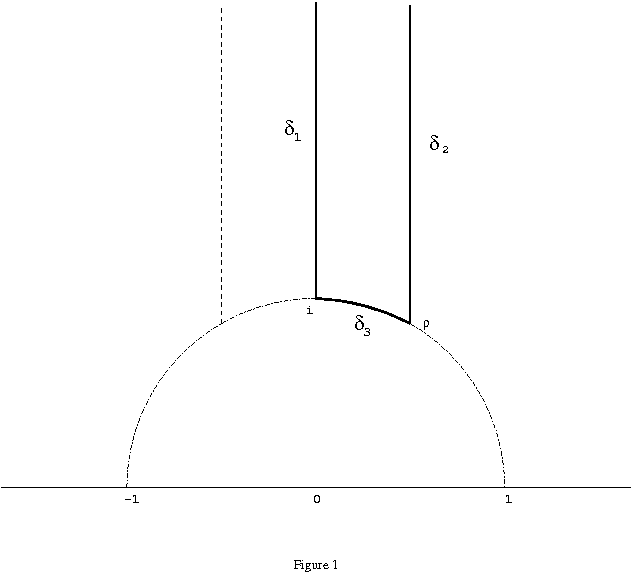}}
\vspace{10pt}
One of the striking consequences of the recent proof of the holomorphic QUE conjecture by Holowinsky and Soundararajan \cite{HS}  is that $\mathcal{Z}(f)$ is equidistributed in $\mathfrak{X}$ as $k \rightarrow \infty$ (see also \cite{Rudnick} and the report \cite{Sa} for a discussion). That is for any nice set $\Omega \subset \mathfrak{X}$
\begin{equation}
\frac{|\mathcal{Z}(f)\cap \Omega|}{|\mathcal{Z}(f)|} \rightarrow \frac{Area(\Omega)}{Area(\mathfrak{X})},
\end{equation}
as $k \rightarrow \infty$. Here $Area$ is the hyperbolic area with $dA =\frac{dx dy}{y^2}$.

This paper is concerned with the zeros of $f$ lying on the geodesic segments $\delta_{1}$, $\delta_{2}$ and $\delta_{3}$ in Figure 1, which we call the `real' zeros of $f$.\footnote{As Zeev Rudnick notes, these segments are the points at which the $j$-invariant is real.}  The reason for this name is that the $a_{f}(n)$'s in (2) are all real and hence $f$ is a real-valued function on the segments $\delta_{1}$ and $\delta_{2}$, while   $z^{\frac{k}{2}}f(z)$ is real-valued on $\delta_{3}$. These follow from the relations
\[
f(S_{1}z) = f(S_{2}z) = \overline{f(z)},
\]
and
\begin{equation}
f(S_{3}z) = \overline{z}^{k}\overline{f(z)}
\end{equation}
where $S_{1}$, $S_{2}$ and $S_{3}$ are the reflections
\[
S_{1}(z)= -\overline{z},\hspace{20pt} S_{2}(z)= 1 -\overline{z},\hspace{20pt} S_3(z)=\frac{1}{\overline z}.
\]
It follows that $\mathcal{Z}(f)$ is invariant under these involutions whose fixed points are the segments $\delta_{1}$,  $\delta_{2}$ and $\delta_{3}$ respectively. Let $\delta^{*} = \delta_{1}\cup \delta_{2} \cup \delta_{3}$. One might expect that the number of real zeros, $N_{real}(f) := |\mathcal{Z}(f)\cap \delta^{*}|$ to grow with $k$, much like the number of real zeros of a random polynomial with real coefficients. According to (3), $N_{real}(f) = o(|\mathcal{Z}(f)|) = o(k)$ as $k \rightarrow \infty$.

\vspace{20pt}
\centerline{\includegraphics[width=3.5in]{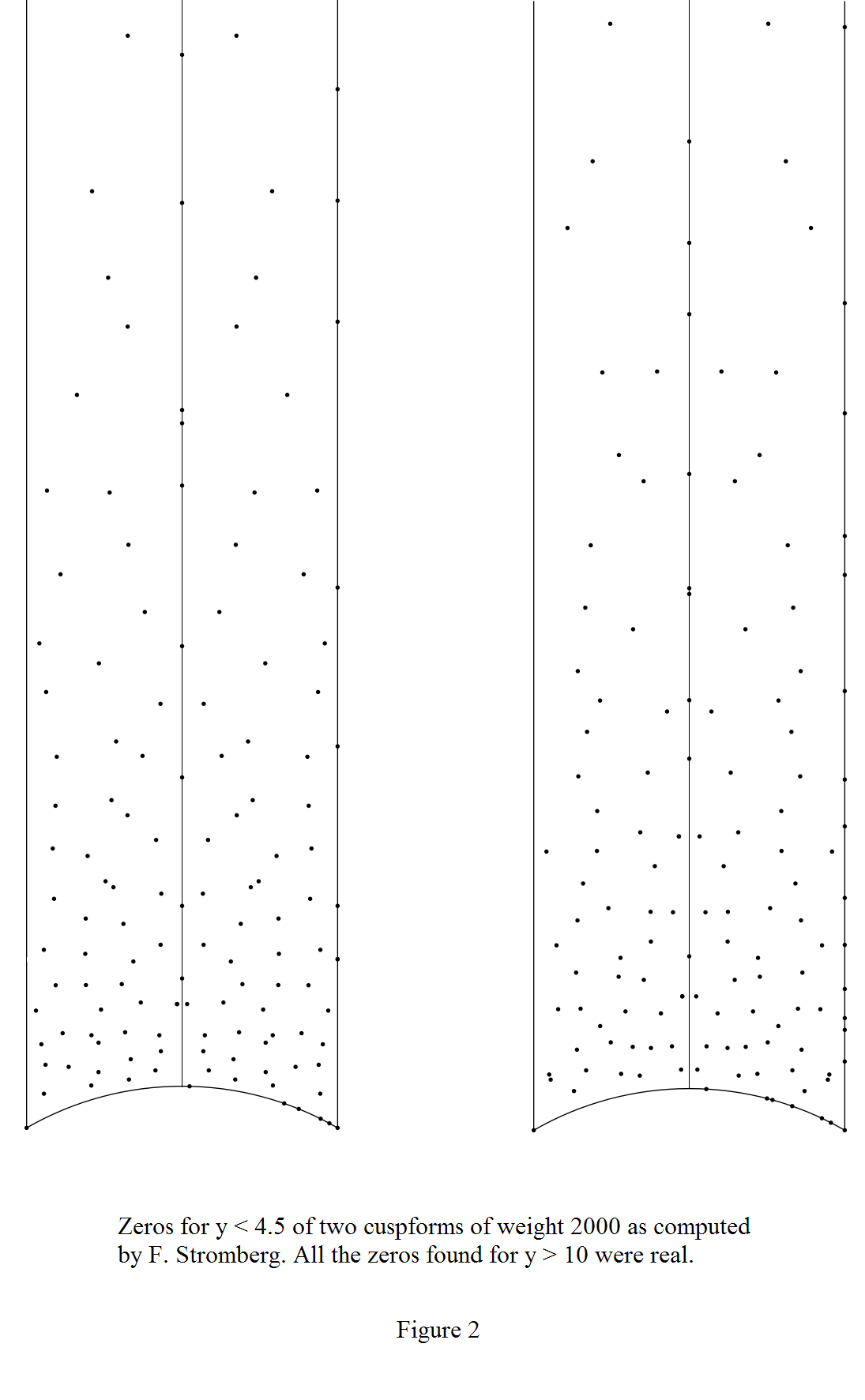}}

\vspace{20pt}

Our first result concerning real zeros of $f$ is that their number does in fact grow with $k$.

\begin{thm}\label{ThmOne} For any $\epsilon > 0$, as $k \rightarrow \infty$ and any Hecke cusp form of even weight $k$,
\[
N_{real}(f) \gg_{\epsilon} k^{(\frac{1}{4} - \frac{1}{80} - \epsilon)}.
\]
\end{thm}

As discussed below and in Section 6, the true order of magnitude of $N_{real}(f)$ is probably $\sqrt{k}\log k$.
\vspace{20pt}

Theorem \ref{ThmOne} follows from a more detailed investigation of the zeros of $f$ in regions which move into the cusp with $k$. To quantify this statement, we define the Siegel sets
\begin{equation}
 \mathcal{F}_{Y} = \{ z \in \mathfrak{X}: \Im (z) \geq Y\},
 \end{equation}
 which have hyperbolic area $\frac{1}{Y}$ and we restrict our considerations to $Y$  satisfying
\begin{equation}
 \sqrt{k\log k} \ll Y < \frac{1}{100}k.
\end{equation}
 We show that the equidistribution (3) holds (approximately) for the zeros of $f$ in these shrinking (relative to the area) regions $\mathcal{F}_{Y}$ and moreover that many of these are real zeros.

\begin{thm}\label{ThmTwo} Let $\epsilon > 0$ and $k \rightarrow \infty$, $f$ any Hecke cusp form of even weight $k$ and $Y$ satisfying (6). Then

\begin{itemize}
\item[(i)]  $$ \frac{k}{Y} \ll |\mathcal{Z}(f) \cap \mathcal{F}_{Y}| \ll \frac{k}{Y}, $$             
\item[(ii)] $$|\mathcal{Z}(f) \cap \mathcal{F}_{Y} \cap (\delta_{1}\cup \delta_{2}| \gg _\epsilon \big{(}\frac{k}{Y}\big{)}^{\frac{1}{2} - \frac{1}{40} -\epsilon}.$$ 
\end{itemize}
The implied constants in (i) above are absolute while that in (ii) depend only on $\epsilon$.
\end{thm}

\vspace{20pt}
Our proof of Theorem \ref{ThmTwo} leads to questions of the nonvanishing of many of the $a_{f}(n)$'s with $n \leq \sqrt{k}$ and it suggests, somewhat surprisingly, that
\begin{equation}
|\mathcal{Z}(f) \cap \mathcal{F}_{Y} \cap (\delta_{1}\cup \delta_{2}| \sim |\mathcal{Z}(f) \cap \mathcal{F}_{Y}|
\end{equation}
as $k \rightarrow \infty$ (see Corollary \ref{Corl 5}). In other words, almost all the zeros in $\mathcal{F}_{Y}$, for $Y$ restricted in the range (6), are real zeros and futhermore half of these are on $\delta_{1}$ and the other half on $\delta_{2}$.

The proof of (ii) in Theorem \ref{ThmTwo} (and so Theorem \ref{ThmOne}) does not specify on which of $\delta_{1}$ or $\delta_{2}$ the zeros that are being produced lie.  On $\delta_{2}$, we are able to take advantage of a natural oscillation induced on the fourier coefficients to produce some zeros of $f(z)$. The situation on $\delta_{1}$ is substantially harder. Our analysis reduces the problem of finding zeros to producing $n$'s with $1 \leq n \ll k^{\frac{1}{2} - \eta}$ and $a_{f}(n) < -\epsilon_{0}$ for some $\eta > 0$ and $\epsilon_{0} >0$ (both fixed independent of $k$). There has been recent progress on the problem of estimating from above the least $n$ for which $a_{f}(n) < 0$, see (\cite{IKS}, \cite{KLSW} and most recently \cite{Mato}).  Remarkably the optimization in \cite{Mato} of the smooth number argument from \cite{KLSW}, together with the sharp subconvex bounds of \cite{Peng} and \cite{J-M} for the critical values of the $L$-function $L(s,f)$, allow us by the closest of margins to produce the requisite $n$'s. In either case (see Section 4), we then have 


\begin{thm}\label{ThmThree} The number of zeros of $f(z)$ on $\delta_{1}$ and separately $\delta_{2}$ goes to infinity as k goes to infinity. More quantitatively,  for $j =1$ and $2$
\[
|\mathcal{Z}(f) \cap \delta_{j}| \gg \log k .
\]
\end{thm}

\vspace{20pt}
It is natural to ask if $N_{real}(f)$ has an asymptotic law. To try to answer this, we determined the number of real zeros such an $f$ would have if the coefficients $a_{f}(n)$ were to behave in some random manner (satisfying the requisite bounds). This model is not completely accurate since the $a_{f}(n)$'s have a multiplicative structure which we ignore (for the sake of simplicity); however we still believe it yields the correct order of magnitude. Whether the constants obtained in the asymptotics are reliable is best checked by numerical experimentation. In any case, such a random model predicts that as $k \rightarrow \infty$
\[
|\mathcal{Z}(f) \cap \delta_{1}| \sim |\mathcal{Z}(f) \cap \delta_{2}| \sim \frac{\sqrt{k}}{4\pi} \log k,
\]
while
\begin{equation}
|\mathcal{Z}(f) \cap \delta_{3}| \sim \frac{\sqrt{k}}{4\pi} \log 3
\end{equation}
and in particular that
\[
N_{real}(f) \sim \frac{\sqrt{k}}{2\pi} \log k.
\]
The random model can also be examined for zeros of $f(z)$ on $\delta_{1}$ and $\delta_{2}$ with $y \gg \sqrt{k}$ (see (26) in Section 6). It predicts that almost all of the zeros of $f(z)$ for $y \gg \sqrt{k}$ are real, which is consistent with the statements made in (7) that were obtained by purely arithmetic considerations. In fact, it is even possible that all of these zeros are real (see remark 5.3). This lends support to the believe that the random model is appropriate even for $y \ll \sqrt{k}$.

\vspace{20pt}
To end this introduction, we outline our proofs of the Theorems. The analysis is based on a suitable approximation to $f(z)$ in the regions $\mathcal{F}_{Y}$ when $k$ and $Y$ are large. This is derived in Sections 2 and 3. Recall that for  a Hecke cuspform, one can write
\[
a_{f}(n) = a_{f}(1)\lambda_{f}(n)n^\frac{k-1}{2},
\]
with $ a_{f}(1)$ nonzero (normalised to $1$) and with $\lambda_{f}(n)$ multiplicative and real.  Among the various inputs into this asymptotic analysis are Deligne's \cite{Deligne} bounds $|\lambda_{f}(n)| \leq d(n) \ll n^{\epsilon}$ for any $\epsilon > 0$. The upshot of the analysis is that for integers $1 \ll l \ll \sqrt{\frac{k}{\log k}}$, with $y_{l} = \frac{k-1}{4\pi l}$ (or for $y$ close enough to $y_{l}$) and $0 \leq \alpha \leq \frac{1}{2}$, both $f(\alpha + iy_{l})$ and $\frac{f'}{f}(\alpha + iy_{l})$ can be approximated by simple functions as long as the $\lambda_{f}(n)$'s are not too small (see Cor 3.3). This condition on $\lambda_{f}(n)$ allows us to conclude that $f$ has exactly $l$ zeros in  $\mathcal{F}_{\frac{k-1}{4\pi l}}$. Thus, part (i) of Theorem \ref{ThmTwo} is reduced to finding some $l$'s in suitable ranges with $\lambda_{f}(l)$ not small. This is a nontrivial problem  since $f$ is varying and $l \ll \sqrt{k}$ , which is small in terms of the conductors of the associated $L$-functions $L(s,f)$ and $L(s,sym^{2}f)$. We proceed by using a much exploited and robust feature, that since $\lambda_{f}(n)$ is multiplicative, either $|\lambda_{f}(p)|$ or $|\lambda_{f}(p^2)|$ is at least $\frac{1}{2}$ for any prime $p$. In section 4 we use this together with a combinatorial analysis to construct sufficiently many such $l$'s. For part (ii) of Theorem \ref{ThmTwo} we proceed by looking for many sign-changes of $f$ on $\delta_{1} \cup \delta_{2}$, and for this we have to elaborate the analysis above by constructing many pairs of $l$'s with opposite parity with $\lambda_{f}(l)$'s not small. We  reduce the problem to seeking a full density set of integers $m$ in $[M,2M]$, for which the short intervals $[m,m+\Delta]$ have at least one prime number. Assuming the Riemann Hypothesis for the Riemann zeta function, Selberg \cite{Sel} showed that the above is true it $\Delta = (\log M)^{2 + \epsilon}$. It appears that the smallest $\Delta$ for which the above is known unconditionally is $\Delta = M^{\frac{1}{20}}$ \cite{Jia}, and this is what we use and it is responsible for the various exponents in our theorems.  As noted earlier, the proof of Theorem \ref{ThmThree} for $\delta_{1}$ relies on strong subconvex bounds for $L(s,f)$, as well as optimised smooth number arguments.

An alternate but equivalent approach is to find $l < l'$ of opposite parity which are close to each other and for which both $\lambda_{f}(l)$ and $\lambda_{f}(l')$ are not too small. This ensures that in the region $\frac{k-1}{4\pi l'} \leq y \leq \frac{k-1}{4\pi l}$, the number of zeros of $f$ is odd and hence by the symmetry associated with $S_{1}$ and $S_{2}$, there must be at least one zero of $f$ on $\delta_{1} \cup \delta_{2}$, in this region.

The last section is devoted to modelling $\lambda_{f}(n)$ by random numbers, namely identical, independently distributed Gaussians of mean zero and variance one. For a random such $f$, we determine the expected density of zeros on each segment of $\delta^{*}$. 

\vskip 0.2in
{\small {\bf Acknowledgements.} 

We thank Fredrik Str\"omberg whose revealing computations and pictures such as those shown in Figure 2, led us to investigate the real zeros of the forms $f$. Thanks also to Andre Reznikov with whom we are preparing a followup to this paper which investigates the `real nodal domains' of Maass forms on $\mathfrak{X}$, to K. Soundararajan for pointing us to the recent preprint \cite{Mato} and to  K. Matomaki for the reference \cite{Jia}.

The first author also thanks the Institute for Advanced Study for providing the possibility 
of an extended visit during which most of this work took place. He also gratefully acknowledges the support from the Ellentuck Fund of the Institute for Advanced Study.}
\vskip 0.2in

\section{Basic proposition.}

We begin with a detailed steepest descent analysis of the behavior of $f(z)$ when $k$ and $y$ are large. In connection with $L^{\infty}$-norms, related approximations are derived in \cite{Sarnak2}(pages 26-29) for Maass forms and by \cite{Xia} for holomorphic forms.
 
Let
\[
I_{s}(y) = y^{\frac{s-1}{2}}e^{-y}
\]
for $y>0$ and $s \in \mathbb{C}$, and define
\[
\Phi_{f}(s;\alpha,y) = \sum_{1}^{\infty} \lambda_{f}(n)e(n\alpha)I_{s}(2\pi ny)
\]
for any real $\alpha$. We then have
\begin{equation}
f(\alpha +iy) = a_{f}(1)(2\pi y)^{-k'}\Phi_{f}(k;\alpha,y),
\end{equation}
where we will use the notation $k' = \frac{k-1}{2}$. More generally, for any $m \geq 1$, we have
\begin{equation}
\big{(}\frac{1}{2\pi i}\big{)}^{m}f^{(m)}(\alpha +iy) = a_{f}(1)(2\pi y)^{-k' -m}\Phi_{f}(k + 2m;\alpha,y).
\end{equation}
This implies that if $f(\alpha +iy) \neq 0$, then
\begin{equation}
\frac{1}{2\pi i} \frac{f'}{f}(\alpha +iy) = \frac{1}{2\pi y}\frac{\Phi_{f}(k+2;\alpha,y)}{\Phi_{f}(k;\alpha,y)}.
\end{equation}
We will first prove our basic 

\begin{prop}\label{Prop 1} Let $\delta > 0$. Then there is a $N(\delta)$ sufficiently large such that for all real $s > N(\delta)$, for all $y$ satisfying $ \sqrt{s} \ll y < \frac{1}{100}s$, and with $B = \sqrt{\delta s\log{s}}$, we have
\begin{equation}
\frac{\Phi_{f}(s;\alpha,y)}{I_{s}(s')} = \sum_{\substack{n \\ |2\pi ny - s'|\leq B}} \lambda_{f}(n)e(n\alpha)e^{-\frac{|2\pi ny -s'|^{2}}{2s'}} + O(s^{-\delta}),
\end{equation}
where $s' = \frac{s-1}{2}$.
\end{prop}

To prove the theorem, we will need some elementary lemmas regarding the behaviour of $I_{s}(y)$.
\begin{lem}\label{Lemma 1.} For a fixed $s$, $I_{s}(y)$ is strictly increasing for $0<y<s'$, and strictly decreasing if $y>s'$.
\end{lem}

\begin{lem}\label{Lemma 2.} Suppose $|h| \ll s^{\frac{2}{3} -\delta}$ for some positive $\delta$ sufficiently small. Then
\[
I_{s}(s' +h) = I_{s}(s')e^{-\frac{h^2}{2s'}}(1 + O(s^{-3\delta})).
\]
\end{lem}
\proof We write
\[
I_{s}(s' +h) = e^{-s' -h}(s')^{s'}(1+\frac{h}{s'})^{s'}.
\]
To a first approximation,
\[
\log(1+\frac{h}{s'})^{s'} = h - \frac{h^2}{2s'} + O(\frac{h^3}{s'^2})
\]
from which the lemma follows.

\vskip 0.1in
{\noindent \bf Proof of Prop. \ref{Prop 1} .} We write, for a paramater $B$ to be chosen later 
\[
\Phi_{f}(s;\alpha,y) = \sum_{i=1}^{3} \Phi_{f}^{(i)}(s;\alpha,y)
\]
where
\[
\Phi_{f}^{(1)}(s;\alpha,y) = \sum_{\substack{n \geq 1\\ 2\pi ny < s' - B}} \lambda_{f}(n)e(n\alpha)I_{s}(2\pi ny),
\]
\[
\Phi_{f}^{(2)}(s;\alpha,y) = \sum_{\substack{n \\ |2\pi ny - s'|\leq B}} \lambda_{f}(n)e(n\alpha)I_{s}(2\pi ny),
\]
and
\[
\Phi_{f}^{(3)}(s;\alpha,y) = \sum_{ 2\pi ny > s' + B} \lambda_{f}(n)e(n\alpha)I_{s}(2\pi ny).
\]
We choose $h= 2\pi ny -s'$ (in Lemma \ref{Lemma 2.}), so that $I_{s}(2\pi ny)=I_{s}(s' + h)$ and consequently assume that $1 \leq B \ll s^{\frac{2}{3}-\delta}$.

We will first estimate $\Phi_{f}^{(3)}(s;\alpha,y)$. Using the fact that $I_{s_{1}}(t) = t^{\frac{s_{1} -s}{2}}I_{s}(t)$, and the upper-bound for $\lambda_{f}(n)$, we have for any $\epsilon >0$ and $s$ sufficiently large,
\[
\Phi_{f}^{(3)}(s;\alpha,y) \ll y^{-\epsilon} \sum_{n > \frac{s' + B}{2\pi y}} I_{s + 2\epsilon}(2\pi ny).
\]
The maximum for $I_{s + 2\epsilon}(t)$ is attained at $s' + \epsilon$ and since $B \geq 1$, we see that 
$I_{s + 2\epsilon}(2\pi ny)$ is decreasing in this sum. We may the approximate the sum by the appropriate integral   
to get 
\[
\Phi_{f}^{(3)}(s;\alpha,y) \ll y^{-\epsilon} \Big{(} \int_{\frac{s' + B}{2\pi y}}^{\infty} (2\pi yt)^{s' + \epsilon}e^{-2\pi yt} \, dt + I_{s + 2\epsilon}(s' + B)\Big{)}
\]
\begin{equation}
\hspace{40pt} \ll y^{-\epsilon} \Big{(} \frac{1}{y}\Gamma(s'+\epsilon +1,s' +B) + I_{s + 2\epsilon}(s' + B)\Big{)}
\end{equation}
where
\[
\Gamma(s,x) = \int_{x}^{\infty} t^{s-1} e^{-t} \, dt
\]
is the incomplete gamma function. 

First, we have
\[
\frac{I_{s + 2\epsilon}(s' + B)}{I_{s}(s')} \ll e^{-B}(1+\frac{B}{s'})^{s'}s^{\epsilon},
\]
so that on using
\[
\log(e^{-B}(1+\frac{B}{s'})^{s'}) = -\frac{B^2}{2s'} +O(\frac{B^3}{s'^2}),
\]
we conclude that 
\begin{equation}
I_{s + 2\epsilon}(s' + B) \ll I_{s}(s')s^{\epsilon}e^{-\frac{B^2}{2s'}}.
\end{equation}
To estimate the incomplete gamma function, we use the following inequality due to Natalini-Palumbo \cite{NP}

\begin{lem}\label{Lemma 3} If $a>1$, $\sigma > 1$ and $x > \frac{\sigma}{\sigma -1}(a-1)$, one has
\[
x^{a-1}e^{-x} < |\Gamma(a,x)| < \sigma x^{a-1}e^{-x}.
\]
\end{lem}
We shall use this lemma with $a=s' +\epsilon +1$ and $\sigma = 1+\epsilon + \frac{s'}{B}$. Then
\[
\Gamma(s'+\epsilon +1,s' +B) \ll \frac{s'}{B}(s' +B)^{s' +\epsilon}e^{-s'-B}
\]
\begin{equation}
\ll \frac{s'}{B}(1+\frac{B}{s'})^{s'}e^{-B}I_{s}(s')s^{\epsilon}
\ll \frac{s'}{B}e^{-\frac{B^2}{2s'}}I_{s}(s')s^{\epsilon}.
\end{equation}
Collecting the estimates from (7),(6) and (5) gives us
\[
\Phi_{f}^{(3)}(s;\alpha,y) \ll y^{-\epsilon}(\frac{s'}{By}+1)e^{-\frac{B^2}{2s'}}I_{s}(s')s^{\epsilon}.
\]

We choose
\[
\sqrt{\delta s\log s} \leq B \ll s^{\frac{2}{3} - \delta}
\]
and $y \gg \sqrt{s}$ to conclude that
\[
\Phi_{f}^{(3)}(s;\alpha,y) \ll I_{s}(s')s^{-\frac{\delta}{2}}.
\]

We next estimate $\Phi_{f}^{(1)}(s;\alpha,y)$. 
Since $n$ is bounded by $s$, we replace $\lambda_{f}(n)$ with $s^{\epsilon}$. Moreover, $I_{s}(2\pi ny)$ is strictly increasing in our interval so that we may approximate the modified sum with the appropriate integral to get
\begin{equation}
\Phi_{f}^{(1)}(s;\alpha,y) \ll s^{\epsilon}\big{(}\int_{1}^{\frac{s'-B}{2\pi y}}(2\pi yt)^{s'}e^{-2\pi yt} \, dt + I_{s}(2\pi y) + I_{s}(s' -B)\big{)}.
\end{equation}
Now, we have 
\[
\frac{I_{s}(2\pi y)}{I_{s}(s')} \ll (\frac{2\pi ye}{s'})^{s'}e^{-2\pi y},
\]
which decays exponentially in $s$ if we choose $y < \frac{s}{100}$. Using the analysis for (6), we have
\[
\frac{I_{s}(s'-B)}{I_{s}(s')} \ll e^{-\frac{B^2}{2s'}}.
\]

To estimate the integral in (8), we break it up into two pieces: let $B_{1} > B$ so that we may write the integral as 
\[
\int_{1}^{\frac{s'-B}{2\pi y}} + \int_{\frac{s'-B_{1}}{2\pi y}}^{\frac{s'-B}{2\pi y}} (2\pi yt)^{s'}e^{-2\pi yt} \, dt.
\]
Replacing $t$ with $\frac{s't}{2\pi y}$ and simplifying gives us 
\begin{equation}
\frac{s'}{y}I_{s}(s')\big{(}\int_{\frac{2\pi y}{s'}}^{1-\frac{B_{1}}{s'}} + \int_{1-\frac{B_{1}}{s'}}^{1-\frac{B}{s'}} t^{s'}e^{s'(1-t)} \, dt \big{)}.
\end{equation}
The integrand is strictly increasing so that the second integral is bounded by 
\[
\frac{B_{1} -B}{s'}(1-\frac{B}{s'})^{s'}e^{B} \ll \frac{B_{1} -B}{s'}e^{-\frac{B^2}{2s'}}.
\]
The first integral is trivially bounded by $e^{-\frac{B_{1}^2}{2s'}}$ provided $B_{1} \ll s^{\frac{2}{3} -\epsilon}$. Hence, from (9), the integral in (8) is 
\begin{equation}
\ll \frac{s'}{y}I_{s}(s')\big{(} e^{-\frac{B_{1}^2}{2s'}} + \frac{B_{1} -B}{s'}e^{-\frac{B^2}{2s'}} \big{)}.
\end{equation}
We choose 
\[
B = \sqrt{\delta s\log s},\hspace{20pt}  B_{1} = \sqrt{\frac{1}{\delta} s\log s}
\]
with $\delta$ sufficiently small to get that (10) is 
\[
\ll \frac{s'}{y}I_{s}(s')(s^{-\frac{1}{\delta}} + \frac{\sqrt{s\log s}}{s'}s^{-\delta}) \ll \frac{\sqrt{s\log s}}{y}I_{s}(s')s^{-\delta} \ll I_{s}(s')s^{-\frac{1}{2}\delta},
\]
since $y>\sqrt{s}$. Collecting all thses estimates together gives us
\[
\Phi_{f}^{(1)}(s;\alpha,y) \ll I_{s}(s')s^{-\frac{1}{3}\delta}.
\]

Finally, for $\Phi_{f}^{(2)}(s;\alpha,y)$ we use Lemma \ref{Lemma 2.}. The error-term contributes
\[
\ll s^{\epsilon}(\sum_{|2\pi ny -s'|\leq B}1)s^{-3\delta} \ll (1+\frac{B}{y})s^{-2\delta} \ll s^{-\delta}.
\] The main-term in Lemma \ref{Lemma 2.} gives the sum stated in (4), with $h=2\pi ny -s'$. This completes our proof.

\section{Main approximation theorem}
Let $l \in \mathbb{N}$ and put $y_{l}(s)= \frac{s-1}{4\pi l}$. Then, 
\[
|2\pi ny_{l}(s) -s'|\leq B \iff |n-l|\leq \frac{Bl}{s'}
\]
so that we must have $n=l$ if $l<\frac{s'}{B}$. Putting in the restrictions on $B$ and $y$ in Prop. \ref{Prop 1} gives us

\begin{thm}\label{Theorem 2} There are positive constants $\beta_{1}$ and $\beta_{2}$ such that for all integers $l$ satisfying $\beta_{1} < l < \beta_{2}\sqrt{\frac{s}{\log s}}$, for all $s$ sufficiently large, the numbers $y_{l}(s)= \frac{s-1}{4\pi l}$ satisfy the equation
\[
\frac{\Phi_{f}(s;\alpha,y_{l}(s))}{I_{s}(\frac{s-1}{2})} = \lambda_{f}(l)e(\alpha l) + O(s^{-\delta})
\]
for some $\delta > 0$, uniformly for any real number $\alpha$.
\end{thm}

We now extend this theorem as follows. Let $m \in \mathbb{N}$. We analyse the behaviour of $\Phi_{f}(s +2m;\alpha,y_{l}(s))$. We may use Prop. \ref{Prop 1} without modification provided $m = o(s)$ so that in the sum in (4), we have
\[
|2\pi ny_{l}(s) -\frac{s+2m-1}{2}| \leq B \iff |n-l+ \frac{ml}{s'}| \leq \frac{Bl}{s'}.
\]
If  we choose $l \ll \epsilon \sqrt{\frac{s}{\log s}}$ and $m \ll \sqrt{s\log s}$, we must have $n=l$ in our sum, so that under these conditions, we have
\[
\frac{\Phi_{f}(s +2m;\alpha,y_{l}(s))}{I_{s+2m}(s' +m)} = \lambda_{f}(l)e(l\alpha)e^{-\frac{m^2}{s+2m-1}} + O(s^{-\delta}).
\]
The exponential term above is $ 1 + O(s^{-\delta})$, provided we choose $m \ll s^{\frac{1}{2} - \delta}$. Next, we see that 
\[
\frac{I_{s+2m}(s'+m)}{I_{s}(s')} = (s' +m)^{m} (1 + O(s^{-\delta})) = s'^{m} +O(s^{m-\delta}).
\]
Combining these estimates gives us
\begin{thm}\label{Theorem3} Let $\delta > 0$ be sufficiently small. For all  $s$ sufficiently large (depending on $\delta$), let $m$ be a real number such that $0 \leq m \ll s^{\frac{1}{2} - \delta}$, and let $l \in \mathbb{N}$ satisfy $ 1 \ll l < \delta \sqrt{\frac{s}{\log s}}$. Then,
\[
\frac{\Phi_{f}(s+2m;\alpha,y_{l}(s))}{I_{s}(\frac{s-1}{2})} = (\frac{s-1}{2})^{m}\lambda_{f}(l)e(\alpha l) + O(s^{m-\delta})
\]
\end{thm}

We apply this theorem with $m=1$ and prove the following
\begin{cor}\label{Theorem 4} Let $\delta > 0$ be sufficiently small. For all  $k$ sufficiently large (depending on $\delta$), let $l \in \mathbb{N}$ satisfy $ 1 \ll l < \delta \sqrt{\frac{k}{\log k}}$ such that $|\lambda_{f}(l)| \gg k^{-\frac{\delta}{2}}$. Then, we have for $y_{l}= \frac{k-1}{4\pi l}$
\[
\frac{1}{2\pi i}\frac{f'}{f}(\alpha + iy_{l}) = l + O(lk^{-\delta}),
\]
uniformly for all $\alpha \in \mathbb{R}$.
\end{cor}
\proof This is a direct consequence of (3) and the theorem above.

\begin{rem} It is easy to see that the theorems above hold with $y_{l}(s)$ replaced by $y$ with $|y-y_{l}(s)| \ll \frac{s^{\frac{1}{2}- 2\delta}}{l}$.
\end{rem}


\section{Sign-changes of \mbox{\boldmath$\Phi_{f}(k;\alpha,y)$}: lower bounds.} 

 We first note that there are no zeros of $f(z)$ if $y > Ck$ for some absolute constant $C$ (one may take $C = \frac{\log 4}{4\pi}$ for example) except for the isolated zero at infinity. This may be deduced directly from the fourier expansion of $f(z)$ by isolating the first fourier coefficient.

To obtain a lower bound for the number of zeros of $f(z)$ on $\delta^{*}$, it sufficies to detect sign-changes of $\Phi_{f}(k;\alpha,y)$ which we recall is  real valued when $\alpha = 0$ or $\frac{1}{2}$ (corresponding to 
$z=\alpha + iy$ lying on $\delta_{1}$ or $\delta_{2}$ respectively). Our theorems in the previous section are valid only for 
\begin{equation}
 \sqrt{k\log k} \ll y < \frac{1}{100}k
 \end{equation}
 and consequently, we restrict our attention to this region. We let $Y$ be a parameter satisfying (11) and we define the Siegel set
  \[
 \mathcal{F}_{Y} = \{ z = \alpha + iy : -\frac{1}{2}<\alpha \leq \frac{1}{2}, y \geq Y\},
 \]
 which is a part of the standard fundamental domain containing the cusp. Then, for $z \in  \mathcal{F}_{Y} \cap \delta^{*}$, we will determine a lower bound for the number of sign-changes of $\Phi_{f}(k;\alpha,y)$ by detecting sign-changes of $\lambda_{f}(l)e(\alpha l)$ and utilising Theorem \ref{Theorem 2}. To this end, we have to ensure that $\lambda_{f}(l)$ is also not too small. 

For the latter, we use the Hecke relations
\begin{equation}
\lambda_{f}(p)^{2} = \lambda_{f}(p^2) +1
\end{equation}
valid for all prime numbers $p \geq 2$. Put $\beta = \frac{\sqrt{5}-1}{2}$. It follows that either $|\lambda_{f}(p)| \geq \beta$ or if not, then $|\lambda_{f}(p^2)| \geq \beta$. We define
\[
\omega = \left\{ 
\begin{array}{ll}
  2 & \quad \mbox{if $|\lambda_{f}(2)| \geq \beta$};\\
  4 & \quad \mbox{if $|\lambda_{f}(2)| <  \beta$},\\ \end{array} \right. 
\]
so that $|\lambda_{f}(\omega)| \geq \beta$.

To obtain a result that is as strong as possible,it would be preferable if we could detect sign-changes of $\lambda_{f}(l)$ for $l$ in suitably short intervals, but as discussed in the introduction, the current methods do not give such sharp results. To circumvent this problem, we look for integers $u$ and $v$ of opposite parity, both in the same short interval such that both $\lambda_{f}(u)$ and $\lambda_{f}(v)$ are not small. The parity condition will then ensure that $\Phi_{f}(k;\alpha,y)$ changes sign for either $\alpha = 0$ or $\alpha = \frac{1}{2}$, for $y$ in a short interval. We will first prove a much weaker result which however has the benefit of localising the detection of sign-changes to each of the lines $\delta_{1}$ and $\delta_{2}$.


\vspace{10pt}
\subsection{Sign-changes on \mbox{\boldmath$\delta_{2}$}.}\

\vspace{10pt}
In this section, we show that there are infinitely many zeros of $f(z)$ on $\delta_{2}$. We will need


\begin{lem}\label{Lemma A} Let $p \geq 2$ be a prime number. For a fixed number $J \geq 1$, there is a constant $B$ depending at most on $J$ and a number $b = b(f,p)$, with $1 \leq b \leq B$ such that if $a = p^b$ we have
\[
\lambda_{f}(a^j) \geq \frac{1}{10}
\]
for all $1 \leq j \leq J$.
\end{lem}
\proof Recall that there is a number $\theta_{p} = \theta(f,p)$ with $0 \leq \theta_{p} \leq \pi$ such that for all non-negative integers $n$,
\[
\lambda_{f}(p^{n}) = \frac{\sin{\big{(}(n+1)\theta_{p} \big{)}}}{\sin \theta_{p}}.
\]
If $\theta_{p} = 0$ or $\pi$, then we may take $ b = 2$ since $\lambda_{f}(p^{2j}) \geq 3$ for all $j \geq 1$. By continuity, we see that $b=2$ still suffices for $\theta_{p}$ near $0$ or $\pi$. In other words, there is number $\theta_{0} > 0$ depending at most on $J$, such that the conclusion of the lemma holds unless $\theta_{0} < \theta_{p} < \pi - \theta_{0}$, which we now assume. 
By Dirichlet's approximation theorem, for any integer $B \geq 1$, there is an integer $1 \leq b \leq B$ such that $\|b\frac{\theta_{p}}{2\pi}\| \leq \frac{1}{B+1}$. Hence, we can find integers $b$ and $b'$ so that $b\theta_{p} = 2\pi b' + \eta$ with $|\eta| \leq \frac{2\pi}{B+1}$ so that
\[
\lambda_{f}(p^{bj}) = \frac{\sin{\big{(}j\eta + \theta_{p} \big{)}}}{\sin \theta_{p}}.
\]
We shall choose $B$ sufficiently large so that $0 < j\eta + \theta_{p} < \pi$ for all $1 \leq j \leq J$. Using estimates for trignometric functions, we conclude that
\[
\lambda_{f}(p^{bj}) \geq 1 - \frac{1}{2}(j\eta)^{2} - j\eta\cot(\theta_{p}) \geq \frac{1}{10}
\]
by choosing $B$ sufficiently large.

\begin{thm}\label{Theorem A} There is a constant $C>0$ such that $f(z)$ has at least $C\log{k}$ zeros on the line $\delta_{2}$ with $z=\frac{1}{2} + iy $ and $ y \geq \sqrt{k\log k}$ for $k$ sufficiently large.
\end{thm}

\proof Let $X \rightarrow \infty$ with $k$ such that $X \ll (\frac{k}{\log k})^{\frac{1}{4}}$, where the implied constant is chosen suitably so as to satisfy the conditions of Theorem \ref{Theorem 2}. We decompose the interval $[1,X]$ into dyadic subintervals $[(2a)^{i},(2a)^{i+1}]$ with $0 \leq i \leq R$ with $R \gg \log k$ and $a \geq 1$ some integer. Each such subinterval we denote by $\mathcal{I}=[m,2am]$ and by $\mathcal{I}^{2}$ the corresponding subinterval $[m^{2},(2am)^{2}]$. Every interval $[m,2m]$ contains a prime number $q \geq 3$, so that both $q$ and $aq$ lie in $\mathcal{I}$. We call a prime $q$ ``good'' if $|\lambda_{f}(q)|\geq \beta$ (see (12)) and ``bad'' otherwise, in which case $|\lambda_{f}(q^{2})|\geq \beta$.

Suppose $q$ is ``good''. In Lemma \ref{Lemma A}, we take $p=2$ and $J=2$ and choose $a = 2^b$ in the dyadic subdivision above. Then $|\lambda_{f}(q)|\geq \beta$, $|\lambda_{f}(aq)|\geq \frac{\beta}{10}$ and
\[
(-1)^{q}\lambda_{f}(q)(-1)^{aq}\lambda_{f}(aq) < 0.
\]
This shows by Theorem \ref{Theorem 2} that $\Phi_{f}(k,\frac{1}{2},y)$ has a sign-change between $\frac{k-1}{4\pi q}$ and $\frac{k-1}{4\pi aq}$.

Now suppose $q$ is ``bad''. In this case, both $q^{2}$ and $a^{2}q^{2}$ lie in $\mathcal{I}^{2}$, $|\lambda_{f}(q^2)|\geq \beta$, $|\lambda_{f}(a^{2}q^{2})|\geq \frac{\beta}{10}$ and 
\[
(-1)^{q^{2}}\lambda_{f}(q^{2})(-1)^{a^{2}q^{2}}\lambda_{f}(a^{2}q^{2}) < 0.
\]
This time we have a sign-change between $\frac{k-1}{4\pi q^{2}}$ and $\frac{k-1}{4\pi a^{2}q^{2}}$. 

By considering only subintervals with $\frac{R}{2} \leq i \leq R$, we can ensure that all our subintervals of the type $\mathcal{I}$ and $\mathcal{I}^{2}$ are disjoint, so that there are at least $\frac{R}{2}$ zeros of $f(z)$ with $z=\frac{1}{2} + iy$ and $y \gg \frac{k}{X^{2}}$, from which the theorem follows.

\vspace{10pt}
\subsection{Sign-changes on \mbox{\boldmath$\delta_{1}$}.}\

\vspace{10pt}
We show in this section

\begin{thm}\label{Theorem20} There is a constant $C>0$ such that $f(z)$ has at least $C\log{k}$ zeros on the line $\delta_{1}$ with $z=iy $ and $ y \geq \sqrt{k\log k}$ for $k$ sufficiently large.
\end{thm}

As mentioned in the introduction, the proof is more involved and we require the following proposition that relies on strong subconvexity estimates for $L$-functions. Our notation here will follow that of \cite{KLSW} and \cite{Mato}.

\begin{prop}\label{Prop 3} There is $\epsilon_{0} > 0$ (independent of $k$) such that if $k$ is large enough and $f$ is of even weight $k$, then there is an $n < k^{0.4963}$ ($n$ a power of a prime) such that 
\[
\lambda_{f}(n) \leq -\epsilon_{0} .
\]

\end{prop}

\proof This follows by a modification of the recent developments connected with the first sign-change in $\lambda_{f}(n)$'s for such $f$'s. Fortuitously, the optimization in \cite{Mato} of the smooth number argument in \cite{KLSW} coupled with the subconvex bounds in \cite{Peng} and \cite{J-M} just allows us to secure an exponent less than $\frac{1}{2}$ in the Proposition.

In more detail, if for $\epsilon > 0$, $\lambda_{f}(p^{e})\geq -\epsilon$ for all prime powers $p^{e} \leq y$, that is $\frac{\sin{\big{(}(m+1)\theta_{p} \big{)}}}{\sin \theta_{p}} \geq -\epsilon$, for $m \leq K$ ($K$ will be fixed, say at 100), then 
\[
\lambda_{f}(p) \geq 2\cos\big{(}\frac{\pi}{m+1}\big{)} - \eta(\epsilon)
\]
for $p \leq y^{\frac{1}{m}}$. Here, we may choose $\eta(\epsilon) = C_{K}\epsilon$ for a constant $C_{K}$ depending at most on $K$ and with $\epsilon >0$ chosen suitable small in what follows. 

Define the multiplicative function $h_{y}$ on squarefree numbers by
\begin{equation}
h_{y}(p) =
\begin{cases}
-2 \hspace{10pt} $if$ \hspace{4pt} p>y, \\
2\cos\big{(}\frac{\pi}{m+1}\big{)} - \eta(\epsilon) \hspace{4pt} $if$ \hspace{4pt} y^{\frac{1}{m+1}} \leq p \leq y^{\frac{1}{m}}, 1 \leq m \leq K\\
2\cos\big{(}\frac{\pi}{K+1}\big{)} - \eta(\epsilon),  \hspace{4pt} p \leq y^{\frac{1}{K+1}}.
\end{cases}
\end{equation}

Following the analysis in \cite{Mato} and in particular the continuity of the solution $\sigma(u)$ to the corresponding difference-differential equation, we conclude that uniformly for $\frac{1}{2} \leq u \leq 3$ and if $\epsilon = \epsilon_{0}$ is small enough but fixed
\begin{equation}
\sum_{n\leq y^{u}} h_{y}(n) \geq \big{(} \sigma_{\epsilon_{0}}(u) + o(1)\big{)}y^{u}
\end{equation}
where $\sigma_{\epsilon_{0}}(u) > 0$ for $\frac{1}{4} \leq u \leq 1.3434 := \kappa$.

As in \cite{KLSW}, define the multiplicative function $g_{y}$ on squarefree numbers by
\begin{equation}
\lambda_{f} = g_{y}*h_{y},
\end{equation}
so that 
\begin{equation}
g_{y}(p) = \lambda_{f}(p) - h_{y}(p).
\end{equation}
Then, by construction $g_{y}(p) \geq 0$ for all p and hence $g_{y}(n) \geq 0$. Now
\begin{equation}
\sideset{}{^\flat }\sum_{n\leq y^{\kappa}} \lambda_{f}(n) = \sideset{}{^\flat }\sum_{d\leq y^{\kappa}}g_{y}(d) \Big{(}\sideset{}{^\flat }\sum_{l\leq \frac{y^{\kappa}}{d}}h_{y}(l)\Big{)},
\end{equation}
where the sums are over squarefree numbers. Since $h_{y}(l) \geq 0$ for $l \leq y^{\frac{1}{3}}$, it follows that $\sideset{}{^\flat }\sum_{l\leq \xi}h_{y}(l) \geq 0$ if $\xi \leq y^{\frac{1}{3}}$ while the same is true for $\xi > y^{\frac{1}{3}}$ from (22). Hence the coefficients of the sum over $d$ in (25) are all non-negative so that
\begin{equation}
\sideset{}{^\flat }\sum_{n\leq y^{\kappa}} \lambda_{f}(n) \geq \sideset{}{^\flat }\sum_{l\leq y^{\kappa}}h_{y}(l) \geq \frac{1}{2} \sigma_{\epsilon_{0}}(\kappa)y^{\kappa}.
\end{equation}
On the other hand, it follows directly from the subconvex bounds for $L(s,f)$ of \cite{Peng} and \cite{J-M} that for any $\delta > 0$
\begin{equation}
\sum_{n\leq y^{\kappa}} \lambda_{f}(n) \ll _{\delta} k^{\frac{1}{3} + \delta}y^{\frac{\kappa}{2}}.
\end{equation}
Combining (26) and (27) leads to a contradiction if $y > k^{\frac{2}{3\kappa} + \delta'}$, which is the case if we assumed the Proposition to be false.

\begin{lem}\label{Lemma14} Given $\xi \geq 1000$ and any cusp form $f$, there are six integers $m$ in the interval $(\xi,50\xi)$ which are relatively prime in pairs and for which 
\[
|\lambda_{f}(m)| \geq \frac{1}{10}.
\]
\end{lem}
\proof The interval $(\sqrt{\xi},\sqrt{50\xi})$ contains at least 18 primes $p$ and for each either $|\lambda_{f}(p)| \geq \beta$ or  $|\lambda_{f}(p^{2})| \geq \beta$ or both (see (20)). If six of these have $|\lambda_{f}(p^{2}| \geq \beta$, then we choose our $m$'s to be these $p^{2}$'s. Otherwise, we can find twelve distinct primes $p_{j}$ with $|\lambda_{f}(p_{j})| \geq \beta$. We now take for our $m$'s the six products $p_{1}p_{2}$, $p_{3}p_{4}$, $...$ , $p_{11}p_{12}$.

\begin{lem}\label{lemma15} Given $\xi \geq 1000$ and $f$, there are relatively prime integers $m_{1}$ and $m_{2}$ in the interval $(\xi,2500\xi)$ such that 
\[
\lambda_{f}(m_{j}) \geq \frac{1}{100},\hspace{10pt} j=1,2.
\]
\end{lem}
\proof Consider the interval $(\sqrt{\xi},50\sqrt{\xi})$. By the previous lemma, there are integers $n_{1}, ..., n_{6}$ in the interval that are relatively prime in pairs such that $|\lambda_{f}(n_{j})| \geq \frac{1}{10}$. Of the three numbers $n_{1}, n_{2}$ and $n_{3}$. at least two have the same sign (we assume the first two) so that $\lambda_{f}(n_{1}n_{2}) \geq \frac{1}{100}$, giving us $m_{1}$ and similarly for $m_{2}$ using the remaining three integers.

\vspace{10pt}
{\noindent \bf Proof of Theorem \ref{Theorem20}.} According to Prop. \ref{Prop 3}, there is an integer $\hat{n}=\hat{n}_{f}$ equal to a prime power $p^{e}$ such that $\hat{n} < k^{0.4963}$ and $\lambda_{f}(\hat{n}) \leq -\epsilon_{0}$ for some fixed $\epsilon_{0} > 0$. Let $\mathcal{I} = (\eta,2500\eta)$ be a subinterval of $(k^{0.4963},\sqrt{\frac{k}{\log k}})$. By Lemma \ref{lemma15}, there is a $m_{1} \in \mathcal{I}$ such that $\lambda_{f}(m_{1}) \geq \frac{1}{100}$. Also applying Lemma \ref{lemma15} but now to the interval $(\frac{\eta}{\hat{n}},2500\frac{\eta}{\hat{n}})$, we find two relatively prime integers $v_{1}$ and $v_{2}$ such that $\lambda_{f}(v_{j}) \geq \frac{1}{100}$ for $j = 1$ and $2$. At least one of the $v_{j}$'s is coprime to $\hat{n} = p^{e}$, say $v_{1}$. We set $m_{2} = \hat{n}v_{1}$ so that $m_{2} \in \mathcal{I}$ and $\lambda_{f}(m_{2}) \leq -\frac{\epsilon_{0}}{100}$. Then, using Theorem \ref{Theorem 2} we find that $f(iy)$ has a sign-change for a $y$ between $\frac{k-1}{4\pi m_{1}}$ and $\frac{k-1}{4\pi m_{2}}$. Since there are $C_{1}\log k$ such disjoint subintervals $\mathcal{I}$ for some positive constant $C_{1}$, we complete our proof.
\vspace{10pt}

\subsection{Sign-changes on \mbox{\boldmath$\delta^{*}$}.}\  

\vspace{10pt}
Let  $X \ll (\frac{k}{\log k})^{\frac{1}{4}}$ be as in the proof of Theorem \ref{Theorem A}. We let $\mathcal{I}$ denote the interval $(X,\sqrt{\omega}X)$ and for any interval $I=(a,b)$, we denote the interval $(\frac{1}{\sqrt{\omega}}a,\frac{1}{\sqrt{\omega}}b)$ by $\frac{1}{\sqrt{\omega}}I$. The intervals $\frac{1}{\sqrt{\omega}}\mathcal{I}$ and $\mathcal{I}$ will form a disjoint pair in our considerations. For any positive $H$ such that $\frac{X}{H} \rightarrow \infty$, let $\mathcal{I}_{j}$ denote the interval $\big{(}X+jH,X+(j+1)H\big{)}$ with $j=0,1,...,R$ with $R$ chosen so that $\mathcal{I}_{j} \subset \mathcal{I}$ for all $j$. We will use the following 

\begin{lem}\label{Lemma 4} There is a $H = H(X) >0$ (as above) such that of the corresponding pairs of subintervals $\{\frac{1}{\sqrt{\omega}}\mathcal{I}_{j},\mathcal{I}_{j}\}$, one can find a positive proportion such that each subinterval of the pair contains at least two (odd) prime numbers.
\end{lem}

We first indicate how we use this lemma to prove
\begin{thm}\label{Theorem31} Let $N_{f}^{Y}(\delta^{*})$ denote the number of zeros of $f(z)$ lying on  $\mathcal{F}_{Y} \cap \delta^{*}$ and let $X$ have order of magnitude $\big{(}\frac{k}{Y}\big{)}^{\frac{1}{2}}$. Then, for  all $k$ sufficiently large, $N_{f}^{Y}(\delta^{*}) \gg \frac{X}{H}$.
\end{thm}
\proof Let $\frac{1}{\sqrt{\omega}}\mathcal{I}_{j_{1}}$ and $\mathcal{I}_{j_{1}}$ be a generic such pair satisfying the lemma so that $\frac{1}{\sqrt{\omega}}\mathcal{I}_{j_{1}}$ contains (at least) two primes denoted by $q$ and $q'$ and $\mathcal{I}_{j_{1}}$ contains $p$ and $p'$. We will construct pairs of integers $\{u,v\}$ with $u$ odd and $v$ even such that neither $|\lambda_{f}(u)|$ nor $|\lambda_{f}(v)|$ are too small. We have 4 cases to consider:
\begin{itemize}
\item[(i)] Suppose $|\lambda_{f}(q)|\leq \beta$ or $|\lambda_{f}(q')|\leq \beta$ (we assume q). If both $|\lambda_{f}(p)|$ and $|\lambda_{f}(p')|$ exceed $\beta$, we put 
\[
u=pp', \hspace{20pt} v=\omega q^{2}.
\]
Then, by the multiplicativity of  $\lambda_{f}(n)$ and the Hecke relations (11)
\[
|\lambda_{f}(u)|\geq \beta^{2},\hspace{20pt} |\lambda_{f}(v)|= |\lambda_{f}(\omega)||\lambda_{f}(q^{2})|\geq \beta^{2}.
\]
Moreover, both $u$ and $v$ lie in the subinterval $\mathcal{J}_{j_{1}}$ where we denote the interval $\big{(}(X+jH)^{2},(X+(j+1)H)^{2}\big{)}$ by $\mathcal{J}_{j}$.
\item[(ii)] If $|\lambda_{f}(q)|\leq \beta$ or $|\lambda_{f}(q')|\leq \beta$ (we assume q) and if either $|\lambda_{f}(p)|$ or $|\lambda_{f}(p')|$ does not exceed $\beta$ (we assume the former), then put
\[
u=p^{2}, \hspace{20pt} v=\omega q^{2},
\]
with the same conclusions as above.
\item[(iii)] Now suppose $|\lambda_{f}(q)|\geq \beta$ and $|\lambda_{f}(q')|\geq \beta$ but with $|\lambda_{f}(p)|\leq \beta$ (or $|\lambda_{f}(p')|\leq \beta$). We choose
\[
u=p^{2}, \hspace{20pt} v=\omega qq',
\]
where this time $|\lambda_{f}(v)| \geq \beta^{3}$ with the condition on $u$ as before.
\item[(iv)] Lastly if all $|\lambda_{f}(.)|$ exceed $\beta$ for the four primes, then we take
\[
u=pp', \hspace{20pt} v=\omega qq',
\]
with estimates involving $u$ and $v$ as before.
\end{itemize}
By Lemma \ref{Lemma 4}, we conclude that there are $\gg R$ disjoint subintervals $\mathcal{J}_{j}$ that contain a pair of integers $\{u,v\}$ such that $|\lambda_{f}(u)|,|\lambda_{f}(v)| \geq \beta^{3}$, with $u$ odd and $v$ even. These disjoint subintervals are contained in the interval $\mathcal{J}=(X^{2},4X^{2})$ so that the conditions of Theorem \ref{Theorem 2} are satisfied with $u$ and $v$ taking the value $l$. We put $y_{u}=\frac{k-1}{2\pi u}$ and $y_{v}=\frac{k-1}{2\pi v}$ and observe that our conditions on $X$ and $Y$ ensure that they exceed $Y$. Using Theorem \ref{Theorem 2} with the two values of $\alpha = 0$ and $\frac{1}{2}$, we see that the size and sign of $\Phi_{f}(k;0,y_{u})$ and $\Phi_{f}(k;0,y_{v})$ are determined by the pair $\lambda_{f}(u)$ and $\lambda_{f}(v)$ respectively. On the other hand, the size and sign of $\Phi_{f}(k;\frac{1}{2},y_{u})$ and $\Phi_{f}(k;\frac{1}{2},y_{v})$ are determined by the pair $-\lambda_{f}(u)$ and $\lambda_{f}(v)$ (due to the parity difference). Consequently, without any further input, we can conclude that at least one of these pairs must be of opposite sign, implying that either $\Phi_{f}(k;0,y)$ or $\Phi_{f}(k;\frac{1}{2},y)$ has a zero for some $y$ between $y_{u}$ and $y_{v}$. Thus there are at least $\frac{1}{2}R$ zeros of $f(z)$ on either $\delta_{1}$ or $\delta_{2}$ and since $R \asymp \frac{X}{H}$, the theorem follows.

\vspace{10pt}
{\noindent \bf Proof of Lemma \ref{Lemma 4}.} To verify the lemma, we appeal to the well-known result that for the interval  $(A,2A)$, for $A$ large enough, there is a number $J$ depending on $A$ such that almost all subintervals of length $J$ in  $(A,2A)$ have at least one prime number. It is easy to check (by combining consecutive subintervals and then discarding subintervals amongst the pairs that do not satisfy the conditions of the lemma) that the lemma is satisfied with $H=\frac{1}{4}J$, for example.

If one assumes the Riemann Hypothesis (for the Riemann zeta-function), then Selberg \cite{Sel} showed in 1943 that one may take $J=g(A)(\log A)^{2}$ for any function $g(A)$ tending to infinity with $A$. This was improved by Heath-Brown \cite{HB} to $J=g(A)\log A$ subject to additional assumptions on the vertical distribution of zeros of the zeta-function. The best uncondtional result to date is due to Jia \cite{Jia}, where $J=A^{\frac{1}{20} + \epsilon}$ (this follows a sequence of similar results by Harman, Watt and Li). As a consequence, we have our

\begin{cor}\label{Corl 1} 
\mbox{}
\begin{itemize}
\item[(i)] Let $N_{f}^{Y}(\delta^{*})$ denote the number of zeros of $f(z)$ lying on  $\mathcal{F}_{Y} \cap \delta^{*}$. Then, for any $\epsilon >0$ and all $k$ sufficiently large, $N_{f}^{Y}(\delta^{*}) \gg \big{(}\frac{k}{Y}\big{)}^{\frac{1}{2} - \frac{1}{40} -\epsilon}$. If one assumes the Riemann Hypothesis for the Riemann zeta-function, then $N_{f}^{Y}(\delta^{*}) \gg \big{(}\frac{k}{Y}\big{)}^{\frac{1}{2}-\epsilon}$ .
\item[(ii)] If $N_{f}(\delta^{*})$ denotes the number of zeros of $f(z)$ lying on  $\delta^{*}$, then for any $\epsilon >0$ and all $k$ sufficiently large, $N_{f}(\delta^{*}) \gg k^{\frac{1}{4} - \frac{1}{80} -\epsilon}$. Moreover, on the Riemann Hypothesis for the Riemann zeta-function, one has $N_{f}(\delta^{*}) \gg k^{\frac{1}{4} -\epsilon}$.
\end{itemize}
\end{cor}


\begin{rem} 
One can give an alternative construction of the numbers $u$ and $v$ used in the proof of Theorem \ref{Theorem31}
above using Lemma \ref{Lemma A} which we do as follows: 
\end{rem} 

Let $X$ be as before and we consider the intervals $\mathcal{I}_{j}$ and $\mathcal{J}_{j}$ as in the proof of Theorem \ref{Theorem31} (we drop the subscripts in what follows). We choose a number $a$ satisfying Lemma \ref{Lemma A} with $p=2$ and $J=4$ and consider the six subintervals $\mathcal{I}$, $a \mathcal{I}$, $a^{2}\mathcal{I}$, $\mathcal{J}$, $a^{2} \mathcal{J}$ and $a^{4}\mathcal{J}$. We also choose $H$ so that there are prime numbers in each of the subintervals $\mathcal{I}$, $a \mathcal{I}$ and $a^{2}\mathcal{I}$, denoted by $p_{1}$, $p_{2}$ and $p_{3}$ respectively. We consider ``good'' and ``bad'' primes as in the proof of Theorem \ref{Theorem A}. Clearly at least two of the primes $p_{1}$, $p_{2}$, $p_{3}$ are ``good'' (which we call case I) or at least two are bad (case II).

In case I, we have that either $a \mathcal{I}$ or $a^{2}\mathcal{I}$ contains two elements $u = a^{j}p$ with $j = 1$ or $2$, and $v=p'$ with both $p$ and $p'$ ``good" odd primes. These numbers are of odd parity and satisfy a lower-bounds of the type described in the proof of Theorem \ref{Theorem31} and so can be used to detect sign-changes. Similarly, in case II, $a^{2} \mathcal{J}$ or $a^{4}\mathcal{J}$ contains two elements $u = a^{2j}p^{2}$ with $j = 1$ or $2$, and $v=(p')^{2}$ with both $p$ and $p'$ ``bad'' odd primes and in this case we get a similar conclusion. The rest of the argument follows that given in Theorem \ref{Theorem31}.


\section{Bounds for the zeros of $\bf{f(z)}$ in \mbox{\boldmath$\mathcal{F}_{Y}$}. }

We now prove a conditional result that gives the precise number of zeros of $f(z)$ in some special Siegel sets, from which we obtain  unconditional upper and lower bounds for the general case.

\begin{thm}\label{Theorem 41} Let $\delta>0$ and k be sufficiently large. Suppose there exists an integer $l$ satisfying $ 1 \ll l < \delta \sqrt{\frac{k}{\log k}}$ such that $|\lambda_{f}(l)| \gg k^{-\frac{\delta}{2}}$. Then there are exactly $l$ zeros of $f(z)$ in the region $\{z=\alpha +iy: -\frac{1}{2}<\alpha \leq \frac{1}{2}, y \geq \frac{k-1}{4\pi l}\}$.
\end{thm}
\proof We integrate $\frac{1}{2\pi i}\frac{f'}{f}(z)$ along the boundary of the indicated region to count the number of zeros (one makes the standard indentations to avoid zeros on the vertical paths). We observe that the vertical integrals cancel due to periodicity and opposite orientations, and that there are no zeros on the horizontal path. Consequently, the number of zeros in our region is
\[
\frac{1}{2\pi i}\int_{-\frac{1}{2}}^{\frac{1}{2}} \frac{f'}{f}(x + iy_{l}) \, dx  \hspace{10pt} = \frac{1}{2\pi i}\big{(} \log f(\frac{1}{2} + iy_{l}) - \log f(-\frac{1}{2} + iy_{l})\big{)}
\]
\[
 = l + O(k^{-\delta}),
\]
by Theorem \ref{Theorem 2}, from which our result follows.


\begin{cor}\label{Corl 5} For sufficiently large $k$, suppose there exists an integer $1 \ll L_{k} \ll \sqrt{\frac{k}{\log k}}$ such that for all integers $1 \leq l \leq L_{k}$, $|\lambda_{f}(l)| \geq \frac{1}{\Delta_{k}}$, where $\Delta_{k} > 0$ also satisfies the condition $\frac{\log \Delta_{k}}{\log k} \rightarrow 0$. Then, all the zeros of $f(z)$ in the region $\mathcal{F}_{Y}$ are real zeros with $Y = A \frac{k}{ L_{k}}$ for a sufficiently large number $A$.

Moreover, if $L_{k}^{+}$ denotes the number of $l$'s above satisfying $\lambda_{f}(l)\lambda_{f}(l+1) > 0$, then $L_{k}^{+}$ of the real zeros above lie on $\delta_{2}$, $L_{k} -1 - L_{k}^{+}$ lie on $\delta_{1}$ and they are all necessarily simple zeros.
\end{cor}
\begin{rem}
From this we expect that asymptotically half of the zeros lie on each segment. As we noted in the introduction, this together with a standard expectation that the $\lambda_{f}(l)$'s are not zero (see for example \cite{F-J}) suggests that perhaps all the zeros of $f(z)$ in $y \gg \sqrt{k}$ are real. There is some numerical evidence provided by Fredrik Str\"omberg that bears this out.
\end{rem}

\proof This follows directly from the above Theorem by looking at the intersection of the Siegel domains for consecutive values of $l$. There is only one zero in such an intersection, and since the zeros are symmetric with respect with the line $\Re(z) = 0$, it must lie on the boundary. Moreover, if $\lambda_{f}(l)\lambda_{f}(l+1) > 0$, then this zero is located on $\delta_{2}$ by Theorem \ref{Theorem 2} and the conclusions follow.

\begin{cor}\label{Corl 2} Suppose $\sqrt{k\log k} \ll Y < \frac{1}{100}k$. Let $N_{f}^{Y}$ denote the number of zeros of $f(z)$  in the region $\mathcal{F}_{Y}$. Then there are absolute positive constants $C_{1}$ and $C_{2}$ such that $C_{1}\frac{k}{Y} \leq N_{f}^{Y} \leq C_{2}\frac{k}{Y}$.
\end{cor}
\proof  Put $X= \sqrt{\frac{k-1}{4\pi Y}}$. By the proof of Theorem \ref{Theorem31}, we see that there are integers $u$ and $u'$ with $u \in (\frac{1}{2}X^{2},X^{2})$ and $u' \in (X^{2},4X^{2})$ satisfying $|\lambda_{f}(u)|, |\lambda_{f}(u')|\geq \beta^{3}$. Then, $ y_{u'} < Y < y_{u}$, where we recall that $y_{l} = \frac{k-1}{4\pi l}$. Applying Theorem \ref{Theorem 41}, we conclude $N_{f}^{Y} \leq N_{f}^{y_{u'}} = u' \ll \frac{k}{Y}$, and $N_{f}^{Y} \geq N_{f}^{y_{u}} = u \gg \frac{k}{Y}$, from which our corollary follows.


\section{A probabilistic model.}

We determine a probabilistic model to predict the expected number of zeros of $f(z)$ on our curves. To consider the distribution of zeros on $\delta_{1}$ and $\delta_{2}$ we may generalise our problem to  any vertical line segment $\mathcal{L}_{\alpha}$ consisting of points $z=\alpha + iy$ with  $y > c$ for any fixed number $c>0$ and $0 \leq |\alpha| \leq \frac{1}{2}$.  We ask for the number of zeros of $f(z)$ as $k$ becomes unbounded. The analysis in the previous sections focused on the part of the lines with $y > \sqrt{k\log k}$, and in fact Corollary \ref{Corl 2} gives us an upper-bound of at most $\sqrt{k}$. In what follows, we will first focus on the range $1 < y \ll \sqrt{k}$.

Writing $\frac{f(\alpha _ iy)}{a_{f}(1)} = \sum_{1}^{\infty} \lambda_{f}(n)e(n\alpha)n^{\frac{k-1}{2}}e^{-2\pi ny}$, we denote the real part by $R_{f}(\alpha,y)$ and the imaginary part by $I_{f}(\alpha,y)$. Since $\sum_{n=1}^{N} \lambda_{f}(n) = o(N)$ and $\sum_{n=1}^{N} |\lambda_{f}(n)|^{2} \sim c(f)N$ for some constant $c(f) > 0$, we replace the $\lambda_{f}(n)$ with independent standard normal random coefficients with mean zero and variance 1 (the constant c(f) plays no role in the subsequent analysis and so may be absorbed in $f$). Following the general principles as shown in Edelman-Kostlan \cite{EK}, we consider the vectors
\[
\mathbf{v} = \mathbf{v}(\alpha,y) =  \sum_{n=1}^{\infty} \cos({2\pi n\alpha}) n^{\frac{k-1}{2}}e^{-2\pi ny} \mathbf{e}_{n}
\]
and
\[
\mathbf{w} = \mathbf{w}(\alpha,y) =  \sum_{n=1}^{\infty} \sin({2\pi n\alpha}) n^{\frac{k-1}{2}}e^{-2\pi ny} \mathbf{e}_{n},
\]
where $\mathbf{e}_{n} = (...,0,0,1,0,0,...)$ denotes the vector with the $1$ in the nth-coordinate. If $\mathbf{u}$ denotes such a vector function, not identically zero, then the probablity density function for the real zeros of the the associated random wave function is given by
\begin{equation}
 \mathcal{P}(\mathbf{u}) = \mathcal{P}(\mathbf{u},\alpha,y) = \frac{1}{\pi} \sqrt{\frac{<\mathbf{u},\mathbf{u}><\mathbf{u'},\mathbf{u'}> - <\mathbf{u},\mathbf{u'}>^{2}}{<\mathbf{u},\mathbf{u}>^{2}}}
\end{equation}
where $\mathbf{u'}$ is the derivative with respect to $y$ and with the standard inner-product. This is the expected number of real zeros of the associated wave function per unit length at the point $y$. We need to compute both $\mathcal{P}(\mathbf{v},\alpha,y)$ and $\mathcal{P}(\mathbf{w},\alpha,y)$. Let 
\begin{equation}
S(k,\alpha,y) = \sum_{n=1}^{\infty} n^{k-1}\cos({4\pi n\alpha})e^{-4\pi ny}.
\end{equation}
Then 
\[
<\mathbf{v},\mathbf{v}>= \frac{1}{2}(S(k,0,y) + S(k,\alpha,y)),
\]
\[
<\mathbf{v},\mathbf{v}'> = -\pi (S(k+1,0,y) + S(k+1,\alpha,y)),
\]
and
\[
<\mathbf{v}',\mathbf{v}'> = 2\pi ^{2}(S(k+2,0,y) + S(k+2,\alpha,y)).
\]
The analogous formulae for $\mathbf{w}$ are the same except the sum of the $S$ functions are replaced with their difference. To determine the asymptotics of $S(k,\alpha,y)$, we apply the Poisson summation formula to the series $S(k,0,y - i\alpha)$ and take real parts so that we get
\[
S(k,\alpha,y) = \Re \frac{\Gamma (k)}{(2\pi i)^{k}} \sum_{h \in \mathbb{Z}} \frac{1}{(2\alpha + h + 2iy)^{k}}.
\]
We truncate this sum with $|h| < 2y$ and see easily that the tailend is bounded by $\frac{y}{(2y)^{k}k}$. If $2\alpha $ is not zero modulo one (that is if $\alpha \neq 0, \pm \frac{1}{2})$) for a fixed $\alpha $, the contribution from $|h| < 2y$ is bounded by $\frac{y}{(2y)^{k}\sqrt{k}}$ so that we conclude that the sum $S(k,\alpha,y)$ is negligible for $1<y<\epsilon \sqrt{k}$. Thus, the inner-products are all determined asymptotically by $S(*,0,y)$ for all $\alpha$. It is easy to evaluate $S(k,0,y)$ using the elementary fact that if a positive function $g(x)$ is increasing in the interval $0 \leq x < x_{0}$ and decreasing for $x > x_{0}$, then 
\[
|\sum_{n=0}^{\infty} g(n) - \int_{0}^{\infty} g(x) \, dx| \ll g(x_{0}).
\]
We apply this with $x_{0} = \frac{k-1}{4\pi y}$, to get
\[
S(k,0,y) = \frac{\Gamma(k)}{(4\pi y)^{k}} \big{(}1+ O(\frac{y}{\sqrt{k}})\big{)},
\]
where here we have used Stirling's formula to estimate the error-term. To apply this result to our distribution function, we have to assume that $\frac{y}{\sqrt{k}}$ tends to zero so that $c<y<\epsilon \sqrt{k}$, with $\epsilon$ positive and sufficiently small. Then we have 
\begin{equation}
\mathcal{P}(\mathbf{u},\alpha,y) \sim \frac{1}{2\pi} \frac{\sqrt{k}}{y}.
\end{equation}
for all fixed $\alpha \neq 0, \pm \frac{1}{2}$ with $\mathbf{u} = \mathbf{v}$ and $ \mathbf{w}$. If $\alpha = 0, \pm \frac{1}{2}$, then $\mathbf{u} = \mathbf{v}$.

\vspace{10pt}

We next model the distribution of zeros on the segment $\delta_{3}$. A direct approach is unnecessary since one can map the segment $\delta_{3}$ onto the line $\Re(z)=\frac{1}{2}$ by the M\"{o}bius transformation $\sigma(z) = \frac{z}{z+1}$. This sends $\delta_{3}$ to the segment $z= \frac{1}{2} + iy$ with $\frac{1}{2} \leq y \leq \frac{\sqrt{3}}{2}$ and the number of zeros of $f(z)$ is the same in both segments. Thus we may apply the analysis above to obtain the same density function (30), which on integrating gives us (8).

\vspace{10pt}
It remains for us to consider the analog of the above when $y \gg \sqrt{k}$. The situation is quite different than in the case  $y \ll \epsilon \sqrt{k}$ since only a few terms dominate  in the sums $S(k,\alpha,y)$, as can be seen by Theorem \ref{Theorem 2}. We will consider the case with $\alpha = 0$ (the case $\alpha=\frac{1}{2}$ is the same) and we will assume $y \gg k^{\frac{1}{2} + \delta}$ with $\delta$ small. We first rewrite $$<\mathbf{v},\mathbf{v}><\mathbf{v'},\mathbf{v'}> - <\mathbf{v},\mathbf{v'}>^{2}$$ in (28) as
\begin{equation}
4\pi^{2}\sum_{m=2}^{\infty}\sum_{n=1}^{m-1} (mn)^{k-1}(m-n)^{2}e^{-4\pi (m+n)y}.
\end{equation}
The summand, considered as a function of real variables $m$ and $n$ has a maximum at the point $m_{0} = \frac{k + \sqrt{k}}{4\pi y}$ and $n_{0} = \frac{k - \sqrt{k}}{4\pi y}$, and then has exponential decay beyond a small neighbourhood of this point, as $k$ becomes unbounded.  The denominator in (28) also localises in a similar manner so that in (29), a maximum occurs at the real value $n= \frac{k-1}{4\pi y}$. Let us put $t = \frac{k-1}{4\pi y}$ and let $l$ denote the integer closest to $t$ in what follows. We will assume that  $y \gg k^{\frac{1}{2}+\delta}$,  so that necessarily  $1\leq l \ll k^{\frac{1}{2}-\delta}$. We observe then that the closest integer to $m_{0}$ and $n_{0}$ is precisely $l$ and see that only the following integer pairs $(l+1,l), (l+1,l-1)$ and $(l,l-1)$ contribute to the sums in (31) while the other terms are exponentially small. Similarly, the integers $l-1, l$ and $l+1$ contribute to $S(k,0,y)$ in (29). We thus conclude that
\[
\mathcal{P}(\mathbf{v},0,y)^{2}  \sim 4 \frac{(1-\frac{1}{l})^{k-1}e^{4\pi y} + 4(1-\frac{1}{l^{2}})^{k-1} + (1+\frac{1}{l})^{k-1}e^{-4\pi y}}{\Big{(}(1-\frac{1}{l})^{k-1}e^{4\pi y} + 1 + (1+\frac{1}{l})^{k-1}e^{-4\pi y}\Big{)}^{2}}.
\]
Our goal is to compute the integral of $\mathcal{P}(\mathbf{v},0,y)$ between say, $Y_{1}<y<Y_{2}$ to give us the expected number of zeros for $y$ on that range on $\delta_{0}$. Thus, we write
\begin{equation}
N := \int_{Y_1}^{Y_2} \mathcal{P}(\mathbf{v},0,y) \, dy = \sum_{l = L_2}^{L_1} \int_{l-\frac{1}{2}}^{l+\frac{1}{2}} \mathcal{P}(\mathbf{v},0,\frac{k-1}{4\pi t}) \Big{(}\frac{k-1}{4\pi t^2}\Big{)} \, dt,
\end{equation}
where $L_{i} \sim \frac{k-1}{4\pi Y_{i}}$.
To evaluate (32) asymptotically, it simplifies our analysis to assume that $k^{\frac{1}{3} + \delta} \ll l \ll k^{\frac{1}{2} - \delta}$, so that we assume $ k^{\frac{1}{2} + \delta} \ll y \ll k^{\frac{2}{3} - \delta}$. This allows us to use approximations of the type 
\[
(1-\frac{1}{l})^{k-1} \sim e^{-\frac{k-1}{l} - \frac{k-1}{2l^2}}
\]
so that we obtain from (32) after some simplifications
\[
N \sim 4 \sum_{l = L_2}^{L_1}(\frac{k-1}{4\pi l^{2}}) \int_{0}^{\frac{1}{2}} \sqrt{\frac{e^{\sigma(\frac{1}{2} - u)} + 4 + e^{\sigma(\frac{1}{2} + u)}}{(e^{-\sigma u} + e^{\frac{1}{2}\sigma}  + e^{\sigma u})^2}}  \, du,
\]
where we made the change of variable $ t = l + u$ and have written $\sigma = \frac{k-1}{l^{2}}$. Noting that $\sigma \rightarrow \infty $, the integral above is well aproximated by 
\[
\int_{0}^{\frac{1}{2}} \frac{e^{\frac{1}{2}\sigma(\frac{1}{2} + u)}}{e^{-\sigma u} + e^{\frac{1}{2}\sigma}  + e^{\sigma u}} \, du \sim \frac{\pi}{2\sigma},
\]
so that the expected number of zeros for $f(iy)$ for $y \geq Y$ with $ k^{\frac{1}{2} + \delta} \ll Y \ll k^{\frac{2}{3} - \delta}$ is $\frac{1}{2}\frac{k}{4\pi Y}$ (see Corollary \ref{Corl 5}). We expect that a more refined analysis will verify this result for all $Y \gg  k^{\frac{1}{2} + \delta}$. Thus the random coefficient model predicts that for $Y \gg k^{\frac{1}{2} + \delta}$
\begin {equation}
|\mathcal{Z}(f) \cap \mathcal{F}_{Y} \cap \delta_{1}| \sim |\mathcal{Z}(f) \cap \mathcal{F}_{Y} \cap \delta_{2}| \sim \frac{1}{2}|\mathcal{Z}(f) \cap \mathcal{F}_{Y}|.
\end{equation}

\vspace{20pt}
\small

\hskip -0.2in
\begin{minipage}{1.8in}
Amit Ghosh\\
Department of Mathematics\\
Oklahoma State University\\
Stillwater, OK 74078, USA\\
e-mail: {\small \bf ghosh@mail.math.okstate.edu}\\
\\ \\ 
\end{minipage}
\hfill
\begin{minipage}{1.8in}
\noindent Peter Sarnak\\
School of Mathematics\\
Institute for Advanced Study\\
Princeton, New Jersey 08540, USA\\and\\
Department of Mathematics\\Princeton University\\
e-mail: {\small \bf sarnak@math.ias.edu}\\
\end{minipage}
\vfill


\vskip 0.5in
\begin{thebibliography}{9}
\footnotesize
\setlength{\parskip=0.0pt}
\setlength{\lineskip=0.0pt} 
 
\bibitem[De]{Deligne} P. Deligne,   {\itshape La conjecture de Weil I},  Publ. Math. IHES , {\bf 43},  (1974), 273-307 

\bibitem[DJ]{DukeJenkins} W. Duke and P. Jenkins, {\itshape On the zeros and coefficients of certain weakly holomorphic modular forms}, Pure Appl. Math. Q. {\bf 4} (2008), no. 4, 1327-1340.

\bibitem[EK]{EK} A. Edelman and E. Kostlan, {\itshape How many zeros of a random polynomial are real?}, Bull. Amer. Math. Soc. (N.S.) {\bf 32} (1) (1995) 1-37.

\bibitem[FJ]{F-J} J. D. Farmer and F. James, {\itshape The irreducibility of some level 1 Hecke polynomials}, Math. of Comp. {\bf 71} (2002), no. 239, 1263-1270.

\bibitem[Ha]{Hahn} H. Hahn, {\itshape On the zeros of Eisenstein series for genus zero Fuchsian groups}, Proc. Amer. Math. Soc. {\bf 135} (2007), no. 8, 2391-2401.

\bibitem[He]{HB} D. R. Heath-Brown, {\itshape Gaps between primes and the pair correlation of zeros of the zeta-function}, Acta Arith. 41 (1982), 85-99.

\bibitem[HS]{HS} R. Holowinsky and K. Soundararajan, {\itshape Mass equidistribution for Hecke eigenforms}, Annals of Math., {\bf 172}, (2010), No. 2, 1517-1528.

\bibitem[Ji]{Jia} C. Jia, {\itshape Almost all short intervals containing prime numbers}, Acta Arith., {\bf 76} (1996), no.1, 21-84.

\bibitem[JM]{J-M} M. Jutila and Y. Motohashi, {\itshape Uniform bounds for Hecke $L$-functions}, Acta Math. {\bf 195}, (2005), 61-115.

\bibitem[IKS]{IKS} H. Iwaniec, W. Kohnen and J. Sengupta, {\itshape The first negative Hecke eigenvalue}, Int. J. Number Theory {\bf 3} (2007), no. 3, 355-363.

\bibitem[KLSW]{KLSW} E. Kowalski, Y. K. Lau, K. Soundararajan and J. Wu, {\itshape On modular signs}, Math. Proc. Cambridge Philos. Soc., {\bf 149} (2010), (3), 389-411.

\bibitem[Ma]{Mato} K. Matom\"{a}ki, {\itshape On signs of Fourier coefficients of cusp forms}, (2010) preprint at  \newline http://users.utu.fi/ksmato/papers/HeckeSigns.pdf .

\bibitem[NP]{NP} P. Natalini and B. Palumbo, {\itshape Inequalities for the incomplete gamma function}, Math. Inequalities and Applications, {\bf 3}, (2000), no. 1, 69-77.

\bibitem[Pe]{Peng} Z. Peng {\itshape Zeros and central values of automorphic $L$-functions}, Ph.D. Thesis Princeton University (2001).

\bibitem[RS]{RanSD} F. K. C. Rankin and H. P. F. Swinnerton-Dyer, {\itshape  On the zeros of Eisenstein series}, Bull. London Math. Soc., {\bf 2} (1970), 169-170.

\bibitem[Ru]{Rudnick} Z. Rudnick, {\itshape On the asymptotic distribution of zeros of modular forms}, Int. Math. Res. Not., {\bf 34}, (2005), 2059-2074.

\bibitem[Sa-1]{Sarnak2} P. Sarnak, {\itshape Letter to Morawetz on $L^{\infty}$ norms of eigenfunctions}, (2004) \newline www.math.princeton.edu/sarnak/ .

\bibitem[Sa-2]{Sa} P. Sarnak, {\itshape Recent progress on the quantum unique ergodicity conjecture}, Bull. Amer. Math. Soc. {\bf48} (2011), 211-228. 

\bibitem[Sel]{Sel} A. Selberg, {\itshape On the normal density of primes in short intervals, and the difference between consecutive primes}, Arch. Math. Naturvid. {\bf 47}, (1943), 87-105.

\bibitem[Ser]{Ser} J-P. Serre, A course in arithmetic, Translated from the French. Graduate Texts in Mathematics, No. 7. Springer-Verlag, New York-Heidelberg, 1973.

\bibitem[Xi]{Xia} H. Xia {\itshape On $L^{\infty}$ norms of holomorphic cusp forms}, J. of Number Theory {\bf 124} (2007), 325-327.

\end{thebibliography}
\end{document}